\documentstyle{article}

\font\emre=msbm10  % scaled \magstep1

\def\Cc{\mbox{\emre C}}
\def\Pp{\mbox{\emre P}}
\def\Rr{\mbox{\emre R}}

\def\Nn{\mbox{\emre N}}
\def\Aa{\mbox{\emre A}}

\newcommand{\k}{\mbox{\bf k}}
\newcommand{\arf}{\mbox{}^{\ast}}
\newtheorem{definition}{Definition}[subsection]
\newtheorem{question}{Question}[subsection]
\newtheorem{answer}{Answer}[subsection]
\def\motto#1{%
         \hbox to\textwidth{\hfill\parbox[t]{18pc}{%
                \verse\footnotesize #1\par
         \vspace{1\baselineskip}}}}
\begin{document}
\title{Arf Rings and Characters\thanks{Expanded
version of an invited talk given at a
  symposium organized at Istanbul to honour Arf's eighty fifth
  birthday.}
}
\author{Sinan Sert\"{o}z \\
Bilkent University, \\ Department of Mathematics, \\
06533 Ankara \\ email: {\tt sertoz@fen.bilkent.edu.tr}}
\date{}

\maketitle

\begin{center}{\it\small
Dedicated to Ordinarius Professor Cahit Arf \\ on the occasion of his
eighty fifth birthday} \end{center}

\begin{abstract}
Algebraic curve branches can be classified according to their
multiplicity sequences. Arf's solution to this problem using Arf
closures and possible implementations of Henselization are discussed.
\end{abstract}

\motto{To learn some modesty one should study curve theory...}

\section{Introduction}
In 1949 an article by Cahit Arf \cite{arf} appears in the Proceedings of the
London Mathematical Society. In this article Arf solves
the classification problem of singular curve branches based
upon their multiplicity sequence. However, the geometric nature of the
problem is so hidden behind the algebraic ideas and subsequent
constructions
that immediately following Arf's article  an article by
Du~Val \cite{duval2} appears, beginning with the
words: ``Cahit Arf's results being
severely algebraic in form,...''. There Du~Val provides the general
reader  with the
necessary geometric ideas that lie behind the scene. This is most
appropriate since it was Du~Val who had formulated the final form of
the problem in his Istanbul article \cite{duval1}
to which Arf refers quite warmly
in his Proceedings article.

In this work I will attempt to describe both the problem
and the solution, using as few machinery as possible, and mostly in
today's terminology, and analyze the `Arf idea' for future prospects.
The problem and its solution form
one of those rare occasions where the solution supplied by an article
answers the question of another article, and does so
without altering the question
to suit the answer! I will describe this solution and then  discuss some
future prospects. However, I must add that neither Arf nor Du~Val can
be held responsible for the yet unsubstantiated optimism that
surrounds the ideas  I will express in the concluding remarks of
this article. I rely on the dexterity of my students to acquit me in
history for any hopes expressed in the final section.

In 1985, fresh out of the graduate school,
I was hired by Erdo\u{g}an \c{S}uhubi as a research assistant to
T\"{U}B{\.I}TAK's Gebze Research Center. The person to whom I would assist
and collaborate with was none other than Cahit Arf himself.
In this article while trying to convey the joy and excitement that
dominated our discussions on this topic I may inadvertently give away
some trade secrets that I learned from him, for which I apologize
from him beforehand. It is a privilege for me to dedicate this
work, albeit humble, to Professor Cahit Arf on the occasion of his
eighty fifth birthday, with gratitude and respect.

\section{The Setup and the Problem}
\subsection{Heuristic Arguments}
The main problem dominating the area is to
understand the behaviour of a curve at its
singularity. Before we attempt any definition however, we must agree to
choose and fix an underlying field to work with. Let us
call our base field \k.
For drawing pictures and extracting intuition \k =\Rr ~
is appropriate but for
most geometric applications \k =\Cc ~is used. Algebraic geometers
generally chose \k ~to be an algebraically closed field of any
characteristic. Arf's arguments however will work for any field \k ~of
any characteristic.
\subsection{The First Example: The Node}
Let us begin with an example. Consider the curve $C$ defined by the
equation ~$y^2=x^2(x+1)$
in the affine plane $\Aa_{\k}^2$. The curve $C$ is the nodal
curve which has a singularity at the origin where the curve intersects
itself. The usual way of `correcting' this singularity is by changing
the space in which the curve lies and hoping that the curve will
behave better given a more suitable environment.
For this we blow up \Aa$^2$ at the origin and
consider the monoidal transform of the curve in the blow up. This
involves first replacing \Aa$^2$ by the new space
\begin{eqnarray}
B_2=\{ ((x,y),[u_1:u_2])\in \Aa^2\times\Pp^1\; |\; u_1y=u_2x\;\}
\end{eqnarray}
which is easily seen to be smooth.
If $\pi :B_2\rightarrow\Aa^2$ is the projection on the $\Aa^2$ component
then $\pi^{-1}(0,0)$ is isomorphic to \Pp$^1$. It is denoted by $E$ and
is called the exceptional divisor. Note that
$\pi^{-1}(\Aa^2-(0,0))$ is isomorphic to $B_2-E$. In particular
$\pi^{-1}(C-(0,0))$ is isomorphic to $C-(0,0)$.
So we have carried the smooth part $C-(0,0)$ of the singular curve $C$
to a new space $B_2$. Now we look for a substitute for the missing point.
A natural way of doing this is by taking the closure
of $\pi^{-1}(C-(0,0))$ in $B_2$, which we denote by $\tilde{C}$.
The question
now is whether $\tilde{C}$ is smooth or not. For this we consider
local coordinates on $B_2$ and examine the equation of $\tilde{C}$ in these
coordinates. Let $U_i$ be the subset of $B_2$ consisting of the points
with $u_i\neq 0$, $i=1,2$. Note that $\{ U_1, U_2\}$ is an open cover
for $B_2$. Define the local
coordinates of $B_2$ as
\begin{eqnarray}
 X=x,\;\; Y=u_2/u_1 \;\;\; {\rm in }\; U_1, \\
X=y, \;\; Y=u_1/u_2 \;\;\; {\rm in }\; U_2.
\end{eqnarray}
With this notation the exceptional divisor $E$ intersects $U_i$ along
the line $X=0$, for $i=1,2$. The equation of  $\tilde{C}$ becomes
\begin{eqnarray}
Y^2=X+1 &~~& {\rm in  }\;\; U_1 \\
1=XY^3+Y^2 &~~& {\rm in  }\;\; U_2.
\end{eqnarray}
We  see that $\tilde{C}$ is smooth in both of these coordinate
neighbourhoods. Moreover  $\tilde{C}$ intersects $E$ at the points
$(0,\pm 1)$ in both charts. The numbers in the y-component
correspond to the slopes of
the tangent lines to $C$ at the origin in $\Aa^2$.
\subsection{The Second Example: The Cusp}
Next consider another example;
Define $C'$ to be the cusp in $\Aa^2$ given by the equation
$y^2=x^5$. Denote by  $\tilde{C'}$  the closure  of
$\pi^{-1}(C'-(0,0))$ in $B_2$. The equation of $\tilde{C'}$ becomes
\begin{eqnarray}
Y^2=X^3 &~~& {\rm in  }\;\; U_1 \\
1=X^3Y^5 &~~& {\rm in  }\;\; U_2.
\end{eqnarray}
We see that   $\tilde{C'}$ intersects $E$ at
the point $(0,0)$ in $U_1$ and is singular there whereas it is smooth in
$U_2$ and does not intersect $E$ there. Judging from the way the
equation of $C'$ is  transformed under the blow up operation we conclude
that if we apply another blow up operation to $\tilde{C'}$ in $U_1$ at
$(0,0)$ then the curve will be transformed to a smooth curve of the
form $Y^2=X$. The  tangent line to the
curve $C'$ at the origin is horizontal and this is reflected in the
$Y$ component of the point where $\tilde{C'}$ intersects $E$ in $U_1$.
\subsection{The General Case}
The crucial information coded in the singularity seems to surface
at the intersection points of the transformed curve with the exceptional
divisor. If we could work with `one piece of information' at a time, then
after each blow up there would be only one intersection with the
exceptional curve and we would continue our analysis from there
on. For this purpose we restrict our attention to such pieces of
information at each singular point on the curve. When we later make
this concept precise we will call it a branch of the curve.

One significant information about the singular point
is the multiplicity of
that point. To find the multiplicity of a point we first count the
number of points a general line intersects the curve. This number is
also known as the degree of the curve. Then we consider a general line
passing through the singular point and count at how many other points
it intersects the curve. We subtract this number from the degree of
the curve and call it the multiplicity of the singular point. This is
reasonable since this difference must count the contribution of that
singular point. Algebraically speaking,
a plane curve is given by a polynomial of
degree $n$ and the number of intersection points of this curve with a
line corresponds to the number of roots of this polynomial after a
linear substitution is made. The number of roots is equal
to the degree of the polynomial when \k ~is algebraically closed.
The situation is similar in $n$-space.

We can thus associate to each singular point its multiplicity. The
multiplicity of a smooth point is 1 by the above definition.
Assume that $p_0$ is a singular point of a curve and that the blow up
of the curve at this point intersects the exceptional divisor at only
one point and further more assume that the same is true for the
subsequent transforms. This property ad hocly describes a singular curve
branch. We
then obtain a sequence $\{ (p_i, m_i)\}$, where $p_i$ is obtained from
$p_{i-1}$ by blowing up and $m_i$ is the multiplicity of $p_i$.
For ease of notation we can only consider the sequence $\{ m_i\}$
which is called the multiplicity sequence of the branch. We now agree
to call two singular branches equivalent if their multiplicity
sequences are the same. The problem is then to classify all singular
branches up to this equivalence class.

\subsection{Technical Formulation}\label{sec:technical}
We observed in the previous sections that the nodal curve $C$ had two parts
to its singularity at the origin whereas the cuspidal curve $C'$
had only one. How can we
recognize this phenomena by looking at their equations? Clearly we wish the
equation of $C$ to split up as the product of two parts and the
equation of $C'$ remain irreducible. The expression $y^2-x^2(x+1)$
is irreducible in the ring $\k [x,y]$. We may say that this ring is
unnecessarily small since it corresponds to global polynomial
functions on \Aa$^2$, whereas we are interested only in what happens
at the origin. Therefore we can look at $k[x,y]_{(x,y)}$, the
localization of $k[x,y]$ at its maximal ideal $(x,y)$. This ring
represents the regular functions at the origin and should fit to our
geometric purpose of focusing our attention to the origin.
However  $y^2-x^2(x+1)$ is still irreducible in this ring.
This
hints to us that we are probably not working in the right rings.
Each irreducible component of  $y^2-x^2(x+1)$
should be of the form $y=\pm x\sqrt{x+1}$. But $\sqrt{x+1}$ is not an element
of the rings  $\k [x,y]$ and $k[x,y]_{(x,y)}$.
Therefore we must find a ring in which
$\sqrt{x+1}$ exists. Observe however that $x+1$ is the square of
1 when computed modulo the maximal ideal corresponding to the
origin. This suggests that we should look at the Henselization of the
local ring $k[x,y]_{(x,y)}$.  (see \cite{milne, david}
for a discussion of Henselization.)
On the other hand the completion of $k[x,y]_{(x,y)}$ with respect to its
maximal ideal always satisfies Hensel's lemma and it can be used at
this stage. In fact in the formal power series ring $\k[[x,y]]$ we can write
$\sqrt{x+1}=\pm (1+x/2-x^2/8+\cdots )$. Hence the equation
$y^2-x^2(x+1)=0$ splits up as $(y-x(1+x/2-x^2/8+\cdots
))(y+x(1+x/2-x^2/8+\cdots )$.

The irreducibility of the expression $y^2-x^2(x+1)$ corresponds to the
fact that
$\k[x,y]/(y^2-x^2(x+1))$, the ring of polynomial functions on
$\tilde{C}$, is an integral domain. We are interested in what happens
at the origin so we localize with respect to the maximal ideal
corresponding to the origin. We can in fact first localize and then
consider the quotient to obtain the ring
$k[x,y]_{(x,y)}/(y^2-x^2(x+1))$.
This is an integral domain. Completing this ring with respect to
its maximal ideal we obtain a ring with zero divisors since
$y^2-x^2(x+1)$ which corresponds to zero can be split up as in the
above discussion. Note however that the equation $y^2-x^5$
continues to stay irreducible even after completing the relevant
ring.
In fact $x$ is never a square in the
ring $R[x]$ where $R$ is a commutative ring with unity. It can only be
a square if $R$ is a suitably chosen noncommutative ring. In our case
the coefficient rings are always fields so $x$ will never split.
(for a discussion of localizations, quotients and completions see
\cite{zariski, atiyah}.)

We can now give a technical definition for a curve branch. Consider
the prime ideal describing an irreducible curve in $n$ space
with a singularity at
the origin. The ideal it generates inside the formal power series
with $n$ indeterminates may split up into components. Each such
component is a branch of the curve passing through the origin. For a
geometric description see \cite{shafarevich}.

\subsection{Du Val's Formulation}\label{sec:duval}
In a much neglected article \cite{duval1} Du Val summarizes the stage
for the classification problem of singular curve branches and
formulates the question whose answer he claims will lead to a
complete understanding of the situation. It is left to the reader to check
that the description of the problem in this section agrees with the
one given in the previous sections.

Define a curve branch $C$ in $n$-space by the following formal
parameterization;
\begin{eqnarray}
x_1  & = & \phi_1(t) \nonumber \\
x_2  & = & \phi_2(t) \nonumber \\
\vdots   & \vdots & \vdots \label{eq:main} \\
x_n  & = & \phi_n(t) \nonumber
\end{eqnarray}
where each $\phi_i(t)$ is a formal power series in $t$
with coefficients from the field \k. We want the branch to pass through
the origin so we impose the condition  that the constant
term of each $\phi_i(t)$ is
zero. Assume that the order of $\phi_1(t)$ is lowest among the
others. We want to blow up the branch at the origin and
write the parameterization of the transformed branch in the coordinate
chart where it intersects the exceptional divisor. The blow up of
\Aa$^n$ ~can be described as
\begin{eqnarray}
B_n=\{ ((x_1,...,x_n), [a_1:\cdots :a_n])\in \Aa^n\times\Pp^{n-1} \;
| \; x_ia_j=x_ja_i, \; 1\leq i,j \leq n\;\}.
\end{eqnarray}
In $U_1$ the local coordinates can be written as
\begin{eqnarray}
X_1=x_1,\; X_2:=a_2/a_1, ... , X_n=a_n/a_1.
\end{eqnarray}
Since $\phi_1(t)$ was chosen to have  the smallest order, $U_1$ is the
chart in which the transformed branch intersects the exceptional
divisor. The parameterization of the transformed branch is now given by
\begin{eqnarray*}
X_1  & = & \phi_1(t) \\
X_2  & = & \phi_2(t)/ \phi_1(t) \\
\vdots   & \vdots & \vdots \\
X_n  & = & \phi_n(t)/ \phi_1(t).
\end{eqnarray*}
Observe that each $X_i$ is again expressed as a formal power series in
$t$. The multiplicity of a branch given by such a representation is
equal to the order of the lowest order series appearing in the
representation.

If the branch  $C$ has its singularity $p_1$ located at
the origin then  its blow up intersects the
exceptional divisor $E_1$ of the first blow up at a point $p_2$.
Semple in \cite{semple} defines  $p_2$ to
be proximate to $p_1$. The second blow up is
centered at $p_2$ and intersects the exceptional divisor $E_2$ of the
second blow up at $p_3$. The transform of $E_1$ also intersects $E_2$
at $p_{12}$. If $p_3=p_{12}$ then we say that $p_3$ is proximate to
$p_1$ and $p_2$. Otherwise it is proximate only to $p_2$. In general
if a certain $p_{i+j}$, with $i,j>0$ lies
on  $E_i$ or  on any transform of $E_i$  under the subsequent blow ups
then we say that $p_{i+j}$ is proximate to $p_{i}$. Nowadays we use
the term ``infinitely close'' instead of ``proximate''. (see
\cite{semple,hartshorne}).
The sum of the
multiplicities of the points $p_1$, $p_2$,...,$p_m$ is called the m-th
multiplicity sum. Du Val in \cite{duval1} a point $p_i$ as
a leading point if the
number of points to which it is proximate is less than the number of
points to which $p_{i+1}$ is proximate. After this bombardment of
definitions comes the most relevant definition: the multiplicity sum
corresponding to a leading point is called a character of the
curve. In other words if $p_i$ is proximate to less points than
$p_{i+1}$ then the sum of the multiplicities of the points $p_1$,...,
$p_i$ is called a character of the branch. The set of all the
characters of the curve was later called the Arf characters. To the
best of my knowledge the first person who first observed the
significance of these numbers and called them ``characters''
was Du Val in \cite{duval1}
and Arf was the first person who could
explicitly calculate them, in \cite{arf}.

Du Val then defines an algorithm, which he calls the modified
Jacobian algorithm, to calculate the multiplicity sequence when these
characters are known. (The Jacobian algorithm, and also the modified
Jacobian algorithm for that matter, calculates the greatest common
divisor of a given set of positive integers.) Finally at the end of
\cite{duval1} he opens a section with the formidable title of
``Outstanding Questions''  and lists some natural questions related to
characters.
\begin{question}[Du Val]How do you find the characters of
a branch if only the local
parameterization with formal power series is available? \label{q1} \\
(see \ref{ans1} for a complete answer).
\end{question}
\begin{question}[Du Val]
Given a set of
positive integers how do you know that they are the characters of some
branch? \label{q2} \\
(see \ref{ans2} for a complete answer).
\end{question}
\begin{question}[Du Val]
Can you find the smallest dimensional space into which the curve branch
can be projected without changing its multiplicity sequence? \label{q3}
\\ (see \ref{ans3} for a complete answer).
\end{question}
Arf recalls that he objected to the amount of
geometric consideration that was clouding the problem,
when Du Val first gave a talk on this subject at
Istanbul University. It must have been 1945.
He claimed that
there was a very algebraic pattern in the problem which could be
solved if one could forget the great geometrical significance of the
problem. Naturally Du Val asked him to work on this. The next day Arf was
homebound with a severe cold so he decided he might as well
think about this problem. Next week when he returned to
work he had in his pocket, scribbled as usual on small pieces of
paper, his own
ticket to immortality...!

\section{The Solution}
In this section we will describe the  tools that are developed by Arf
to solve the above problem. In the course of these descriptions it may
seem to the reader that we have strayed away from the problem. But
despite mounting evidence against we will be doing geometry and all
this will be justified at the end when we describe how these pieces
fall in to complete the jigsaw.

\subsection{Generalities and Some Notation}\label{sec:general}
We will be working in the formal power series ring $\k[[t]]$
of a single indeterminate $t$. If $H$ is a subring of $\k[[t]]$ then
we define
\begin{eqnarray}
W(H) & = & \{ ord \alpha \; |\; \alpha \in H\; \} \\
     & = & \{ i_0=0<i_1<\cdots <i_r<\dots\;\}
\end{eqnarray}
The integers $i_0,i_1,...$ form a semigroup of the additive group of
nonnegative integers \Nn. We assume that $H$ is always so
chosen  that the
semigroup $W(H)$ contains all integers large enough. In other words if
$\nu_l$ denotes the greatest common divisor of the integers
$i_1,i_2,...,i_l$, then for $\rho$ large enough we want
$\nu_{\rho}=1$. This is not a serious restriction since if $H$ does not
satisfy this condition then $H$ can be transformed by an automorphism
of $\k[[t]]$ into a subring $H'$ which satisfies this condition.
Arf does not carry out this transformation  but chooses
and fixes a appropriate $T\in\k[[t]]$ whose order is $\nu=gcd W(H)$,
and assumes throughout that his ring $H$ can be considered as a subring
of the power series ring in the variable $T$ if necessary, see \cite[p
258, Remarque]{arf}.

For each $i_r$ in
$W(H)$ let $S_{i_r}$ be an element of $H$ with ord$S_{i_r}=i_r$.
We define an ideal $I_h$ by
\begin{eqnarray}
I_h&=&\{\alpha\in H\; |\; ord\alpha \geq h\;\}.
\end{eqnarray}

It can be shown that the inverse of any element in $H$ of order zero
is again an element of $H$. With this in mind we define the {\bf set}
$I_h/S_h$ as
\begin{eqnarray}
I_h/S_h &=& \{ \alpha S_h^{-1}\in H\; |\; \alpha\in I_h\;\}.\label{eq:arf1}
\end{eqnarray}
This set consists of certain elements of $H$ closed under addition
since $I_h$ is an ideal. But the product of any two elements from this
set need not be in the set. We want to consider the smallest subring of $H$
containing the set $I_h/S_h$. It turns out that this ring is
independent of the particular element $S_h$ we have chosen. So we
define the {\bf ring}
\begin{eqnarray}
[I_h]&=&\mbox{\rm the smallest subring of $H$ containing $I_h/S_h$}.
\label{eq:arf2}
\end{eqnarray}
Similarly for a semigroup $G=\{i_0=0<i_1<i_2<\cdots \}$
of nonnegative integers and for $h\in G$ define
\begin{eqnarray}
G_h &=& \{\alpha\in G\;|\; \alpha \geq h\;\} \\
G_h-h &=& \{  (\alpha-h)\in G\;|\; \alpha \geq h\;\} \label{eq:arf3}\\
{[} G_h{]} &=& \mbox{\rm the semigroup of nonnegative integers}
 \nonumber \\
& &\mbox{generated by the set $G_h-h$} .\label{eq:arf4}
\end{eqnarray}
\subsection{Arf Rings and Arf Semigroups}
It is clear now that only for special rings $H$ the set $I_h/S_h$ will
already be a ring. We single out such rings and give them a name:\\
\begin{definition}[Arf Ring] {\rm \cite[p 260]{arf}} \\
{A subring $H$ of the formal power series ring
  $k[[t]]$ is
  called an Arf ring if the set $I_h/S_h$ is a ring for every nonzero
  $S_h\in H$. (i.e.  $I_h/S_h=[I_h$ for every $S_H\in H$. See
 the equations  \ref{eq:arf1} and \ref{eq:arf2}.)}
\end{definition}
Similarly we select out special semigroups: \\
\begin{definition}[Arf Semigroup] {\rm \cite[p 260]{arf} } \\
{ A semigroup $G$ of nonnegative integers is
  called an Arf semigroup if the set $G_h-h$ is a semigroup for every
  $h\in G$. (i.e. $G_h-h=[G_h]$ for every $h\in G$. See the equations
  \ref{eq:arf3} and \ref{eq:arf4}.)
  \ref{eq:arf2})}
\end{definition}
What if a ring $H$ does not satisfy this condition? We then associate
to it a ring which does:
\begin{definition}[Arf Closure of a Ring] {\rm \cite[p 263]{arf} }\\
{If $H$ is a subring of
    $\k[[t]]$ then we define the Arf closure $\arf H$ of $H$ to be the
    smallest Arf ring in $\k[[t]]$ containing $H$.}
\end{definition}
We similarly define Arf closure for semigroups:
\begin{definition}[Arf Closure of a Semigroup] {\rm \cite[p 263]{arf}}\\
{If $G$ is a
    subsemigroup of the nonnegative integers \Nn$=\{0,1,2,...\}$ then
    we define the Arf closure $\arf G$ of $G$ to be the smallest Arf
    semigroup in \Nn ~containing $G$.}
\end{definition}
It remains to check that these definitions are not void. The
ring $\k[[t]]$ itself is obviously an Arf ring. So for any subring $H$
the collection of Arf rings in $\k[[t]]$ containing $H$ is not empty
and by  Zorn's lemma must have a smallest element which we call $\arf
H$. Here ordering is done with respect to inclusion. Similarly the
semigroup \Nn ~is an Arf semigroup and thus the definition of Arf
closure for semigroups is not void.

\subsection{Arf Characters}
In the previous section we saw that we are building parallel
constructions in algebra and arithmetic but it was not clear from
their definitions that they would interact in a meaningful way.
We are now ready to observe a crucial interaction.
If $H$ is a subring of the power series ring, as described in section
\ref{sec:general}, first take its Arf closure $\arf H$ and then
look at $W(\arf H)$, the semigroup of orders of the Arf closure. There
is a smallest semigroup $g_{\chi}$ in \Nn ~whose Arf closure is equal to
$W(\arf H)$. This semigroup $g_{\chi}$ has a minimal generating set
$\chi_1,...\chi_n$ which generates it over \Nn. These are the Arf
characters of $H$:
\begin{definition}[Arf Characters] {\rm \cite[p 265]{arf}}\\
If $H$ is a subring of
  $\k[[t]]$ as described in \ref{sec:general}, then the characteristic
  semigroup of $H$ is defined to be the smallest semigroup $g_{\chi}$ in \Nn
  ~with $\arf g_{\chi}=W(\arf H)$. The semigroup $g_{\chi}$
  can be generated over
  \Nn ~with a minimal set of positive
   integers $\chi_1,...,\chi_n$ which are
  defined to be the characters of $H$.
\end{definition}
As is mentioned in section \ref{sec:duval} the concept of characters
seems to have originated with Du Val's article \cite{duval1}. However
the following ideas appear for the first time in Arf's article
\cite{arf} for the explicit purpose of solving the problem raised in
question \ref{q3}.
\begin{definition}[Bases, Base Characters, Dimension] {\rm \cite[pp
  271-274]{arf} } \label{def:bases}  \\
  If $H$ is an Arf ring let $X_1$ be an
  element in $H$ of smallest positive order. $X_1,...,X_{n-1}$
  having been chosen let $X_n$ be an element of smallest order in $H$
  not included in the Arf closure of the ring
  $\k[X_1,...,X_{n-1}]$. (The ring $\k[X_1,...,X_r]$ is defined as
  consisting of the elements of the form $\sum\alpha_{j_1j_2\cdot
    j_r}X_1^{j_1}X_2^{j_2}\cdots X_n^{j_r}$ where
  $\alpha_{j_1j_2\cdots j_r}\in\k$ and the summation is taken over all
  $(j_1,j_2,...,j_r)\in\Nn^r$.)
  Since $W(H)$ is finitely generated over \Nn
  ~this process terminates and we obtain  a finite collection of elements $\{
  X_1,...,X_m\}$. Such a collection is called a base of $H$. If we
  denote their orders by $\chi_i=ord(X_i)$ for $i=1,...,m$, then the
  numbers $ \chi_1,...,\chi_m$ are called the base characters of
  $H$. The number $m$ is called the dimension of $H$.
\end{definition}
These definitions are backed up  with concrete
constructions. If $G$ is a semigroup of \Nn ~then a method of
constructing its Arf closure is described in \cite[p 263]{arf}.
If $H$ is a subring of $\k[[t]]$ the a method of constructing its Arf
closure is given in \cite[p 267]{arf}. With these methods the
characters can be calculated in a finite number of steps. But not
content with this Arf gives a direct method of computing the
characters of any Arf semigroup, \cite[p 277]{arf}.
\subsection{Solving the Problems}
The setup being as in section \ref{sec:duval} we now explain the
procedure for answering the questions \ref{q1}, \ref{q2} and \ref{q3}.
At our disposal we only have some elements $\phi_1(t),..., \phi_n(t)$
of $\k[[t]]$ coming from the parameterization of the branch $C$ as in
the equation
\ref{eq:main}. We now answer the question \ref{q1}.
\begin{answer} Denote by $H$ the ring $\k[\phi_1(t),..., \phi_n(t)]$
  generated by the $\phi_i(t)$'s in the formal power series ring
  $k[[t]]$. (See the definition \ref{def:bases}
for a description of the ring
  $\k[\phi_1(t),...,\phi_n(t)]$.) First construct its Arf closure $\arf H$.
  and then construct the smallest
  semigroup $g_{\chi}$ whose Arf closure if $W(\arf H)$. The minimal
  generators of $g_{\chi}$ are the characters of the given branch $C$.
  Or alternatively use the method described in \cite[p 277]{arf} to
  find the characters directly from $W(\arf H)$.
\label{ans1}
\end{answer}
The validity of this answer is proved in \cite[p 266, Th\'eor\`eme
3]{arf} where Arf shows that Du Val's Jacobian algorithm applied to
these characters gives the sought for multiplicity sequence.

Now we come to the  question \ref{q2} of knowing
whether a given set of positive integers $ 0<\gamma_1<\cdots <\gamma_l$
are characters of an actual branch.
\begin{answer} If $G$ denotes the semigroup of \Nn ~generated by
  $ 0<\gamma_1<\cdots <\gamma_l$, then the characters of $\arf G$ form a
  subset of $\{ \gamma_1,...,\gamma_l\}$, \cite[p 256, Th\'eor\`eme
  2]{arf}.  If any element from the
  semigroup $\arf G$ is added to this set of characters then the
  resulting set of integers will give the same multiplicity sequence
  when Du Val's modified Jacobian algorithm \cite[p 108]{duval1}
  is applied to them, \cite[p 266, Th\'eor\`eme 3]{arf}.
  Once $\arf G$ is known one can construct
  elements $\phi_1(t),..., \phi_r(t)$ of $\k[[t]]$ such that the ring
  $\arf H=\k[\phi_1(t),...,\phi_l(t)]$ is an Arf ring and that
  $W(\arf H)=\arf G$, \cite[p 277 and p 282, Th\'eor\`eme 6]{arf}.
  \label{ans2}
\end{answer}

And we finally come to the  question \ref{q3} of finding the
smallest dimensional space into which the branch $C$ given by the
equations  \ref{eq:main}
can be projected without changing its multiplicity sequence.
\begin{answer} If $H$ denotes the ring $\k[\phi_1(t),...,
  \phi_n(t)]$, defined as in the definition \ref{def:bases}, then
  the space of smallest dimension into which the
  branch $C$ can be projected without changing its multiplicity
  sequence has its dimension equal to  the dimension of $\arf
  H$. (Recall the definition of the dimension of $\arf H$ as
  given in the definition
  \ref{def:bases}.)
  \label{ans3}
\end{answer}
To prove the validity of this answer Arf defines in \cite[p 273]{arf}
a generating system as a set of elements $Y_1(t),...,Y_{\nu}(t)$ in
$H$ such that the Arf closure of
$\k[Y_1(t)-Y_(0),...,Y_{\nu}(t)-Y_{\nu}(0)]$ is equal to $\arf
H$. Then on \cite[p 274]{arf} shows that the smallest number of elements
required in a generating system is the dimension of $\arf H$. And
finally in \cite[p 279, Th\'eor\`eme 5]{arf} he shows how to construct
a ring $\arf H$ of a given set of characters such that $\arf H$ has
the smallest possible dimension.

For the reader who wants to see only a geometric argument summarizing
all this we refer to the last section of Arf's article
\cite[pp 285-287]{arf} where he
relates all this ``algebraic interpretation'' to the problems raised
by Du Val in \cite{duval1}.

\subsection{Silent Heros: Base Characters}

In the statements of the above answers it was not necessary to quote
even the existence of base characters let alone their importance.
However they play a crucial
role in shedding light into the whole scene besides actually answering
question \ref{q3}. To understand the role they play in this set up
first observe that  if $\arf H$ is an Arf ring then $\arf
H_h:=[I_{i_h}]$
is also an Arf ring. (See section \ref{sec:general} for the
notation. In particular note that $W(\arf H)=\{ i_0<i_1<i_2<\cdots \}$
and hence the subscript $h$ counts the number of possible
constructions up to that stage.)  Arf shows that the characters of
each $\arf H_h$ are determined by characters of $\arf H$ but the base
characters of each $\arf H_h$ constitute a new set of characters. And
he then sets out to demonstrate how they can be constructed, \cite[pp
274-275]{arf}. Thus he obtains for each $\arf H$ a set of invariants:
the characters of $\arf H$ and the base characters of $\arf H_h$ for
each $h\geq 0$.
Moreover in\cite[p 282, Th\'eor\`eme 6]{arf} he shows
that if $l_c$ denotes the number of characters of $\arf H$ and if $l_b$
denotes the least possible number of base characters of $\arf H$
(which he shows
how to construct in \cite[p 279, Th\'eor\`eme 5]{arf}) then for any
integer $n$ with $l_b\leq n\leq l_c$ there exists an Arf ring who has
the same characters and is of dimension $n$. To show that all this is
possible he calculates a concrete example in \cite[pp 283-284]{arf}
where he starts with an Arf ring $\arf G$ and first calculates its
characters. Then he calculates all possible sets of invariants that
can be associated to $\arf G$, as he explained how to do in \cite[pp
277-278]{arf}, and finishes the example by producing an actual $\arf
H$ for each set of invariants having that set as its set of
invariants. So a careful rereading of Arf's article shows that the
algebraic structure of a branch is totally understood with the help of
base characters. Despite their importance they have never been given
the ``Arf'' adjective which they silently deserved...

\section{Concluding Remarks}
The most significant follow up of Arf rings came with Lipman's 1971
article \cite{lipman} in American Journal of Mathematics. In fact it
seems that Lipman was the first mathematician to coin the expression
``Arf Rings'' in the literature. He seems to have been motivated by
the similarity of ideas in Arf rings and in Zariski's theory of
saturations. In this article he studies the condition
that Arf singled out for his rings and relates them to Zariski's ideas in
analyzing the singularities.

The basic idea seems simple; if there is
a singularity then the local ring there `misses' something and the
idea is to fill these `gaps' in a controlled manner so as to
understand the nature of the singularity. A similar idea in a much
grand scale was utilized by Hironaka in \cite{hironaka} where he
measures how far his local rings are from being regular. Then Bennett
showed that in higher dimensions Hilbert functions can be used
effectively to measure these `gaps', see \cite{bennet}. Arf's idea,
being simple and fundamental surfaces every now and then in the study
of singularities. Recently the Italian and Spanish mathematicians are
working on ideas around Arf rings.

As already mentioned in section \ref{sec:technical} the completion of
the local ring at the singularity may be too large. Henselization may
be enough to find a complete set of invariants but the calculations
may be involved, to say the least. In higher dimensions instead of checking
the deviation from regularity by Hilbert functions I believe that a
direct description of the `gaps' may be much more useful. The recent
developments on Gr\"obner bases provides a margin of hope that this is
possible. This particular speculation
aims to provoke several definitions trying
to pinpoint a feature in the ring as a `gap'...
On the other hand the innocent looking structure of a branch involves
questions about the structure of complete rings which are for some
cases answered satisfactorily by Cohen in \cite{cohen}, however the
theory is far from being complete.

When the topic is on curves one can
speculate forever! But past the basic definitions comes a land of no
man... I remember what I heard years ago in a conference on curves:
``To learn some modesty one should study curve theory...''

{\small \bf Acknowledgments:}{\small ~I thank Professors G. Ikeda and
  A. H. Bilge for inviting me to give a talk at the symposium for the
  occasion of Professor Arf's eighty fifth birthday. The occasion has
  rekindled my interest and awakened fond memories. Professor Lipman
  has kindly communicated to me his informal thoughts and memories on the
  topic. Professor Eisenbud has answered a crucial question about
  complete rings, which saved a lot of time and worry for me in the
  preparation process of this manuscript. Even though I freely used
  the information I received from them the same cannot be said about their
  wisdom which I no doubt failed to transfer to the manuscript due to
  my own inability to recognize a wise word when I see one... Last but
  not least my thanks are to Professor Arf for being a motivation and
  a convincing example for us.}

\end{document}